\documentclass[12pt]{article}
\usepackage{amsfonts}
\usepackage{amsmath}
\usepackage{amssymb}
\usepackage{amsthm}
\usepackage{color}
\DeclareMathAlphabet{\mathpzc}{OT1}{pzc}{m}{it}

\newcommand{\ed}{\mathrm{d}}
\newcommand{\const}{{\rm const}}

 \newcommand{\ovl}[1]{\overline{#1}}

\newcommand{\wh}[1]{\widehat{#1}}

\newcommand{\reff}[1]{(\ref{#1})}
\newcommand{\pr}[1]{\vspace{.5cm}\textit{Proof of Theorem} \ref{#1}}
\newcommand{\Proof}{\textit{Proof }}

\numberwithin{equation}{section}
\newtheorem{theorem}{Theorem}

\newtheorem{lemma}{Lemma}

\theoremstyle{remark}
\newtheorem{rem}{Remark}[section]

\theoremstyle{definition}

\newcommand{\R}{{\mathbb R}}

\renewcommand{\P}{{\sf P}}
\newcommand{\E}{{\sf E}}

\newcommand{\eps}{\varepsilon}
\title{Gibbs Random Graphs}
\author{Pablo A. Ferrari,$^{~3}$
\and Eugene A.~Pechersky,$^{~1}$
  \and Valentin V.~Sisko$^{~2}$
  and Anatoly A.~Yambartsev$^{~3}$
   }

\begin{document}
\maketitle {\footnotesize

\noindent $^1$Dobrushin laboratory of Institute for Information
Transmission Problems of Russian Academy of Sciences,\\ 19,
Bolshoj Karetny, Moscow, Russia.\\
E-mail: pech@iitp.ru

\noindent $^2$Federal University of Fluminense, Institute of
Matematics, Rua Mario Santos Braga, s/n,
24020-140 Niter\'oi, RJ, Brazil.\\
E-mail: valentin@mat.uff.br

\noindent $^3$Department of Statistics, Institute of Mathematics
and Statistics, University of S\~ao Paulo, Rua do Mat\~ao 1010,
CEP 05508--090, S\~ao Paulo SP, Brazil.\\
E-mail: yambar@ime.usp.br}

\begin{abstract}
  Consider a discrete locally finite subset $\Gamma$ of $R^d$ and the complete
  graph $(\Gamma,E)$, with vertices $\Gamma$ and edges $E$. We consider Gibbs
  measures on the set of sub-graphs with vertices $\Gamma$ and edges $E'\subset
  E$.  The Gibbs interaction acts between open edges having a vertex in
  common. We study percolation properties of the Gibbs distribution of the graph
  ensemble. The main results concern percolation properties of the open edges
  in two cases: (a) when the $\Gamma$ is a sample from homogeneous Poisson process and (b)
  for a fixed $\Gamma$ with exponential decay of connectivity.



\end{abstract}

\section{Introduction}
Let a sample $\Gamma\subset\R^{d}$ of a \emph{point process} be a
locally finite set of $\R^{d}$. We consider an ensemble of graphs
whose \emph{vertices} are the points of $\Gamma$ and whose
\emph{edges} are the set of unordered pairs of points in $\Gamma$.
Each edge can be open or close; we study probability distributions
on the set of configurations of open edges. The classical example
is the \emph{Erd\"os-R\'enyi's random graph} where each edge is
open independently of the others with some probability. Many works
have studied this model and correlates (see for example \cite{BB},
\cite{RD}). In this paper we introduce interactions between edges
and/or vertices and study the associated Gibbs measures.

Given a configuration of open edges, we say that two edges
\textit{collide} if both of them are open and they have a vertex
in common. We call \emph{monomers} those vertices that are extreme
of no open edge. A positive energy is paid by each collision and
by each monomer. Furthermore to any open edge it is assigned a
positive energy proportional to its length. This energy function
is described explicitly in \eqref{H} below. Roughly speaking,
measures associated to this energy function give more weight to
configurations with few monomers, no collisions and short edges.

In Theorem~\ref{exist} we prove the existence of an infinite
volume Gibbs measure associated to the energy function.  Call
\emph{dimer} an open edge not colliding with any other edge. In
Theorem~\ref{ground} we show that the ground states for a typical
configuration $\Gamma$ is composed only by monomers and dimers.
Theorem~\ref{uniqgr} gives conditions for the uniqueness of the
ground state. Then we consider the problem of percolation for two
distributions of the point process.  In the first case we consider
that $\Gamma$ is a sample of Poisson process and prove that if the
density is small, then there is no percolation
(Theorem~\ref{Poiss}). Second we ask $\Gamma$ to satisfy a
$\varepsilon_0$-hard core condition (that is, the ball of radius
$\varepsilon_0$ around each point of $\Gamma$ has no other point
of $\Gamma$), where $\epsilon_0$ is an arbitrary small positive
number and prove that there is no percolation with probability 1
in this case (Theorem~\ref{core}).

To prove Theorem~\ref{core} we construct a process called \textit{cluster
  branching process} and use a coupling of the paths of the process and the
configurations of the Gibbs random graph.

The open question is the non-percolation for Poisson vertices set
having arbitrary large rate $\lambda$. Of course, the large
density of Poisson vertices must be compensated by a large
temperature.

\section{Definitions }

Let $\Gamma$ be sample of a point process and consider the
complete graph $G=(\Gamma, E)$, where $\Gamma$ is the set of
vertices and $E$ is the set of all unordered pairs
$[\gamma_{1},\gamma_{2}]\subset\Gamma$. The set $E$ can be
represented as $E=\bigcup_{\gamma\in\Gamma}E_{\gamma}$ where
$E_{\gamma}$ is the set of edges $e\in E_{\gamma}$ such that
$\gamma\in e$. The length of edge $e=[\gamma_{1},\gamma_{2}]$ is
defined by $L(e)=|\gamma_{1}-\gamma_{2}|$.

Let $\ovl{\Omega}=\{0,1\}^E $ be the set of configurations $\omega:\
E\rightarrow \{0,1\}$. The edge $e\in E$ is called \textit{open} with respect to
$\omega$ if $\omega(e)=1$, and the edge is \textit{closed} if
$\omega(e)=0$. Often we shall use the term open/close with no mentions the
configuration $\omega$.

We define the graph $G(\omega)=(\Gamma, E(\omega))$ as the
subgraph of $G$ whose edges are $E(\omega)$, the open edges of
$\omega$. The \emph{degree} $d_{\omega}(\gamma)$ of a vertex
$\gamma\in \Gamma$ with the respect of $\omega\in \ovl{\Omega}$ is
$d_{\omega}(\gamma)=\#\{e\in E_{\gamma}:\:\omega(e)=1\}$, the
number of open edges containing $\gamma$. We shall use both the
configuration $\omega$ and the graph $G(\omega)$ as the synonyms.

Our goal is to define a Gibbs distribution on the ensemble
$${\Omega}=\{\omega:\: \omega\in
\ovl{\Omega},\;d_{\omega}(\gamma)<\infty \mbox{ for all }\gamma \}$$ the set of
graphs whose vertices have finite degree.  Introduce the following formal
Hamiltonian
\begin{equation}\label{H}
H(\omega)=\sum_{e:\:\omega(e)=1}L(e)+\sum_{\gamma\in \Gamma}
\phi_{\omega}(\gamma),
\end{equation}
where $\phi_{\omega}(\gamma)$ is a ``penalty'' energy function defined by
\begin{equation}\label{phi}
\phi_{\omega}(\gamma)= \left\{ \begin{array}{ll} h_0,
& \mbox{if }d_{\omega}(\gamma)=0; \\
0,& \mbox{if }d_{\omega}(\gamma)=1; \\ h_1
\binom{d_{\omega}(\gamma)}2 & \mbox{if }d_{\omega}(\gamma)>1.
\end{array}\right.
\end{equation}
where $h_0$ and $h_1$ are fixed positive parameters. Notice that
$\phi_{\omega}(\gamma)$ depends only on the degree $d_{\omega}(\gamma)$. It
defines the energy of a pair interaction between open edges from $E_{\gamma}$.

An aspect which is not standard for Gibbs field constructions is that the
potential function $\phi_{\omega}$ depends on infinite number of `sites'. The
edges play here the role of sites in a lattice. To know that
$\phi_{\omega}(\gamma)=0$ we have to check that $\omega(e)=0$ for infinite many
$e\in E_{\gamma}$. However the generalization of the usual Gibbs construction
for this case is rather easy and does not cause special
considerations. Therefore further we do not concern this peculiarity.

The description of the infinite volume Gibbs measure with
Hamiltonian $H$ requires the definition of finite volume Gibbs
measures. Taking a finite volume $\Lambda\subset\R^{d}$ the set of
points $\Gamma_{\Lambda}=\Gamma\cap\Lambda$ is finite. The
complete finite graph $G_{\Lambda}=(\Gamma_{\Lambda},E^{\Lambda})$
has edges $E^{\Lambda}$ connecting all pairs of points of
$\Gamma_{\Lambda}$. Let $\Omega_{\Lambda}$ be the set of
configurations $\omega:\:E^{\Lambda}\to \{0,1\}$. The Gibbs state
$\P_{\Lambda}$ on $\Omega_{\Lambda}$ with the `free' boundary
condition is defined by
\begin{equation}\label{Statefree}
\P_{\Lambda}(\omega)=\frac{\exp\{-\beta H_{\Lambda}
(\omega)\}}{Z_{\Lambda}},
\end{equation}
for $\omega\in \Omega_{\Lambda}$, where the parameter $\beta$ is the inverse
temperature, $Z_{\Lambda}$ is the normalizing constant and
\begin{equation}\label{Hlambdafree}
H_{\Lambda}(\omega)=\sum_{e\in
E^{\Lambda}(\omega)}L(e)+\sum_{\gamma\in \Gamma_{\Lambda}}
\phi_{\omega}(\gamma),
\end{equation}
where the set $E^{\Lambda}(\omega)\subseteq E^{\Lambda}$ is the set of all open
edges in $\Lambda$.

Since $\ovl{\Omega} $ is compact, there exists a Gibbs distribution on
$\ovl{\Omega} $ which may have infinite-degree vertices.  We show later
finiteness of the degrees with probability 1. This implies that the Hamiltonian
\ref{Hlambdafree} generates an infinite volume Gibbs field concentrated on
$\Omega$.

\section{Main results}
\subsection{Existence} A point set $\Gamma$ is \textit{weakly homogeneous} if
for any $\gamma\in \Gamma$ and any $\beta>0$
\begin{equation}\label{hom}
  T_{\gamma}(\beta)=\sum_{e\in E_{\gamma}}e^{-\beta
    L(e)}<\infty,
\end{equation}
where $E_{\gamma}$ is the set of all possible edges with common
extreme $\gamma$.
A point set $\Gamma$ is \textit{strongly homogeneous} if
\begin{equation}\label{hom1}
 \sup_{\gamma\in\Gamma}T_{\gamma}(\beta)<\infty.
\end{equation}

Since the possible unbounded contribution to the sum \eqref{hom} comes from
accumulation of vertices that are close to $\gamma$, if $\Gamma$ consists of
hard core ball centers of a fixed radius, then $\Gamma$ is strongly
homogeneous. We show later in Lemma \ref{poisfin} that for a Poisson process
with law $\pi_{\lambda}$, almost all $\Gamma$ is weakly homogeneous but not
strongly homogeneous.

\begin{theorem}\label{exist}
  For any weakly homogeneous $\Gamma$ and any $0<h_{0}<h_{1}$ the Gibbs random
  graph distributions associated to the Hamiltonian $H$ defined in \eqref{H} are
  concentrated on $\Omega$.
\end{theorem}
Notice that the theorem concerns all possible Gibbs measures $\P$ associated to
the Hamiltonian $H$; the problem of uniqueness is not discussed in this article.

\subsection{Ground states}
A configuration $\wh{\omega}$ is a local perturbation of $\omega\in\Omega$
if $\wh{\omega}\neq \omega$ and there exists a finite volume
$\Lambda\subset\R^{d}$ such that $\wh{\omega}$ coincides with $\omega$ for edges
not included in $\Lambda$: $\wh{\omega}(e)= \omega(e)$ for all
$e\not\subset\Lambda$.  A configuration $\omega\in\Omega$ is a \textit{ground
  state} if for any local perturbation $\wh{\omega}$ of $\omega$
\begin{equation*}
H(\wh{\omega})-H(\omega)\geq 0.
\end{equation*}
(The difference is well defined because all but a finite number of terms vanish.)

\begin{theorem}\label{ground} For any finite-local $\Gamma$ and any
  $0<h_{0}<h_{1}$ there exists at least one ground state.  Furthermore if
 $\omega$ is a ground state of the Gibbs random graph distribution, then
 $$d_\omega(\gamma) \le 1$$ for every $\gamma \in
\Gamma$. Moreover, the length of edges in a ground state are less
than $2h_0$.
\end{theorem}

\begin{theorem}\label{uniqgr}
  Let $\pi_{\lambda}$ be the distribution of a homogeneous Poisson process
    with rate $\lambda>0$. There exists $\lambda_{g}$ such that if
    $\lambda<\lambda_{g}$, then for $\pi_{\lambda}$-almost all $\Gamma$ the
    ground state of the Gibbs random graph with vertices $\Gamma$ is unique.
\end{theorem}

\subsection{Non-percolation at low rate or low temperature} Let
$\omega\in\Omega$. The set $E(\omega)$ is split into a set of
maximal connected components which we call \textit{clusters}.  For
a locally finite $\Gamma$, we say that the associated to $\Gamma$
Gibbs random graph measure $\P$ on $\Omega$ \textit{percolates} if
there exists an infinite cluster in $E(\omega)$ with $\P$
probability 1. Otherwise we say that the Gibbs random graph
measure $\P$ on $\Omega$ does not percolate.

In the next theorem we establish non percolation of $\P$ when $\Gamma$ is a
Poisson process with small intensity $\lambda$.
\begin{theorem}\label{Poiss}
  Let $\pi_{\lambda}$ be the distribution of a Poisson process with
  rate
  $\lambda>0$ and $\Gamma$ chosen with $\pi_{\lambda}$. Then in the region
\begin{equation}\label{region}
F=\left\{(\lambda,T):\ \lambda\leq\frac{1}{2h_0+J(T)}\right\},
\end{equation}
where
    \[J(T)=\int_{2h_0}^{\infty}\frac{e^{-x/T}}{e^{-x/T}+e^{-2h_0/T}}{\rm d}x,
\]
the Gibbs measure $\P$ associated to $\Gamma$ does not percolate,
$\pi_{\lambda}$-almost surely.
\end{theorem}

The next theorem is stronger but for more restricted sets $\Gamma$.
\begin{theorem}\label{core}
  Let $\Gamma$ be strongly homogeneous. Then there exists a critical temperature
  $T_{c}(\Gamma)$ such that for all $T<T_{c}(\Gamma)$ the Gibbs measure $\P$
  associated to $\Gamma$ does not percolate.
\end{theorem}

\section{Proofs}
\subsection{Main Lemma} Let $\gamma
$ be a point of $\Gamma$, and $\Sigma_{\gamma
}$ be a set of all configurations defined on $E_{\gamma
}$ and having a finite degree at $\gamma
$. That is, any $\sigma\in \Sigma_{\gamma}$ is the restriction of a
configuration $\omega\in \Omega$ to $E_{\gamma
},$ $\sigma=\omega_{E_{\gamma
  }}$; in this case we call $\sigma$ a \textit{star} centered at $\gamma
$, or simply a star. Clearly $\sigma\subset\omega$. Let $d_{\sigma}:=
d_{\omega}(\gamma
)$.

Let $\Lambda\subset\R^d$ be a finite volume and 
$\gamma
\in \Lambda$; let $E_{\gamma
}^{\Lambda}$ be the set of edges contained in $\Lambda$ having $\gamma
$ as its end, and $\Sigma_{\gamma
}^{\Lambda}$ be the set of the configurations on $E_{\gamma
}^{\Lambda}$ which are restrictions of the configurations from
$\Omega_{{\Lambda}}$.

\begin{lemma}\label{main}
Let the point $\gamma\in\Lambda$, consider $\sigma\in\Sigma_{\gamma
}^{\Lambda}$ and
$\Omega_{{\Lambda}}(\sigma)=\{\omega\in\Omega_{{\Lambda}}:
\:\sigma\subset\omega\}
$, and assume
that the number of open edges $|\sigma|$ in $\sigma$ is greater or
equal than 2 then
\begin{equation}\label{estim}
\P_\Lambda(\Omega_{{\Lambda}}(\sigma)) \le e^{-\beta \sum_{e:
\sigma(e)=1 } L(e) - \beta h_1 \binom{d_\sigma} {2} } e^{\beta h_0
(d_\sigma + 1) }
\end{equation}
\end{lemma}
\proof Let $V(\sigma)$ be the set of all vertices belonging to open edges of the
star-configuration $\sigma$
$$ V(\sigma) = \{ \gamma': \ [ \gamma,\gamma' ]
\in E_{\gamma}^{\Lambda}(\sigma)\}.
$$
To any configuration $\omega\in\Omega_{{\Lambda}}(\sigma)$ we associate a
configuration $\tilde\omega$ without the star $\sigma$; in $\tilde\omega$ the
point $\gamma$ is isolated:
\begin{equation}\label{wh}
\tilde{\omega}(e)=\begin{cases} {\omega}(e)& \mbox{ if }e\notin
\sigma,\\ 0, &\mbox{ if }e\in \sigma.
\end{cases}
\end{equation}
The transformation $\omega\mapsto\tilde \omega$ (taking out the star $\sigma$
from $\omega$) changes the penalty weights only at the vertices in
$\{\gamma\}\cup V(\sigma)$. Namely, for any point $\gamma'$, the penalty weight
in configuration $\tilde \omega $ is
\begin{equation} \label{new-pen} \phi_{\tilde{\omega}}(\gamma')=
\begin{cases} \phi_{{\omega}}(\gamma')&\mbox{ if }\gamma'\notin V(\sigma)
\cup \{\gamma\},\\
h_{1}\binom{d_{\omega}(\gamma')-1}2+h_{0}\delta_{0}(d_{\omega}(\gamma')-1)&
\mbox{ if }\gamma'\in V(\sigma), \\ h_0 & \mbox{ if
}\gamma'=\gamma,
\end{cases}
\end{equation}
where $\delta_{0}(\cdot)$ is Kronecker symbol. Consider the
possible changes of the energy caused by the star $\sigma$ removal
from configuration $\omega.$ The difference between energies is
\begin{equation} \label{dif}
H(\tilde\omega) - H(\omega) = - \sum_{e\in \sigma} L(e) +
\sum_{\gamma' \in V(\sigma) \cup \{ \gamma\}}
(\phi_{\tilde\omega} (\gamma') - \phi_\omega (\gamma'))
\end{equation}
and
\begin{eqnarray}
&& \phi_{\tilde\omega} (\gamma') - \phi_{\omega} (\gamma')
\label{delta} =  \\ && \phantom{aaaaaaa} \nonumber
\begin{cases}
h_0 - h_1 \binom{d_\sigma}{2} & \mbox{ if } \gamma'=\gamma, \\
- h_1 (d_{\omega} (\gamma') - 1) & \mbox{ if } \gamma'\in V(\sigma)
\mbox{ and } d_\omega(\gamma') \ge 2, \\
h_0 & \mbox{ if }\gamma'\in
V(\sigma) \mbox{ and } d_\omega(\gamma') =1.
\end{cases}
\end{eqnarray}
We used here the fact that for any $\gamma'\in V(\sigma)$ its
degree in configurations $\omega$ and $\tilde \omega$ satisfy
equality $d_{\tilde\omega} (\gamma') = d_{\omega} (\gamma') -1$, and
we used that $\binom{k}{2} - \binom{k-1}{2} = k-1$ for $k\ge 2.$
Let us denote
 $$\Delta \phi_\omega(\gamma'):= \phi_{\tilde\omega}
(\gamma') - \phi_{\omega} (\gamma'). $$
 Then the probability of the star $\sigma$ is
\begin{eqnarray}
 \P_\Lambda(\Omega_{{\Lambda}}(\sigma)) &=&
\sum_{\omega\in\Omega_{\Lambda}(\sigma)} \P_\Lambda (\omega) =
\frac{1}{Z_{\Lambda}} \sum_{\omega\in\Omega_{\Lambda}(\sigma)}
e^{-\beta H(\omega)} \nonumber \\
&=& \frac{1}{Z_{\Lambda}}
\sum_{\omega\in\Omega_{\Lambda}(\sigma)} e^{-\beta
(H(\omega)-H(\tilde\omega))} e^{-\beta H(\tilde\omega)} \nonumber
\\
&=& \frac{1}{Z_{\Lambda}} \sum_{\omega\in\Omega_{\Lambda}(\sigma)}
e^{ -\beta \sum_{e\in \sigma} L(e) + \beta \sum_{\gamma' \in
V(\sigma) \cup \{\gamma\} } \Delta \phi_\omega (\gamma') }
e^{-\beta H(\tilde\omega)} \label{est1}
\end{eqnarray}
Here and further instead of the sum $\sum_{e: \sigma(e)=1}$ we
write simply $\sum_{e\in \sigma}.$

In the last expression of (\ref{est1}) the factor $e^{-\beta H(\tilde\omega)}$
does not depend on $\sigma$ and the factor $e^{ -\beta \sum_{e\in \sigma} L(e)}$
does not depend on $\omega$.  However the factor $e^{\beta \sum_{\gamma' \in
    V(\sigma) \cup \{\gamma\} } \Delta \phi_\omega (\gamma') }$ depends on both
$\sigma$ and $\omega$. To find an upper bound depending only on
$\sigma$ we represent the energy of the difference as
\begin{eqnarray*}
\sum_{\gamma'\in V(\sigma)\cup \{\gamma\}}
\Delta\phi_{\omega}(\gamma') &=& h_0-h_1\binom{d_\sigma}{2} +\sum_{\begin{smallmatrix}\gamma'\in
V(\sigma)\\
d_\omega(\gamma')=1
\end{smallmatrix}}h_0-\sum_{\begin{smallmatrix}\gamma'\in
V(\sigma)\\
d_\omega(\gamma')>1
\end{smallmatrix}}h_1(d_\omega(\gamma')-1) \\
&\leq &h_0+h_0d_\sigma-h_1\binom{d_\sigma}{2}.
\end{eqnarray*}
The above inequality we obtain if we assume that  all vertices in
$V(\sigma)$ have its degrees equal to 1.

Thus we obtain the following estimate for
$\P_\Lambda(\Omega_{{\Lambda}}(\sigma))$. Let $\Omega(\gamma)$
denote the set of configurations where $\gamma$ is isolate point,
then it follows from \reff{est1} that
\begin{equation}\label{est4}
\P_\Lambda(\Omega_{{\Lambda}}(\sigma)) \le e^{ -\beta \sum_{e\in
\sigma} L(e) + \beta (h_0 - h_1 \binom{d_\sigma}{2}) +\beta h_0
d_\sigma } \frac{1}{Z_{\Lambda}} \sum_{\tilde\omega \in
\Omega(\gamma)} e^{-\beta H(\tilde\omega)}.
\end{equation}
Noting that
 $$ \frac{1}{Z_{\Lambda}} \sum_{\tilde\omega \in
   \Omega(\gamma)} e^{-\beta H(\tilde\omega)} = \P_\Lambda (\Omega(\gamma)) <
 1$$ we obtain the estimation (\ref{estim}) of the lemma.  \qed

\medskip

\begin{rem} The estimate \reff{estim} does not depend on
$\Lambda$ when $V(\sigma)\subset\Lambda$.
\end{rem}

\medskip

We can generalize the lemma for the case when there is an
``environment''. Let $B$ and $F$ be some nonempty sets of edges of
$E^{\Lambda} _{\gamma}$ without intersection $B\cap F=\emptyset$.
Let $B$ be the set of open edges and $F$ be the set of closed
edges. Introduce a configuration of the ``environment'' $\mu$ on
$B\cup F$
\begin{equation*} \mu(e)=\begin{cases} 1, &\mbox{ if } e\in B \\
0,&\mbox{ if } e\in F \end{cases}
\end{equation*}
Consider the following sets of edges $G_{\gamma,\;B\cup
F}^{\Lambda} = E_{\gamma} ^{\Lambda}\setminus (B\cup F).$ And let
$\Sigma_{\gamma,B\cup F}^\Lambda$ be the set of all configurations
on $G_{\gamma, \;B\cup F}^{\Lambda}$. If $\sigma \in
\Sigma_{\gamma,B\cup F}^\Lambda$ then the degree of the point
$\gamma$ is equal to the number of open edges on $\sigma$ (we
denote it by $d_\sigma$) plus the number of the open edges $|B|$
in the ``environment'' $\mu$ (denote it by $d_\mu$). As before
denote $\Omega_\Lambda (\sigma)$ and $\Omega_\Lambda(\mu)$ the
sets of all configurations which include the star-configuration
$\sigma$ and the configuration $\mu$ correspondingly.
\begin{lemma}\label{main1}
Let $\sigma\in \Sigma_{\gamma_{0}, B\cup F}^{\Lambda}$ then
\begin{equation}\label{ineq1}
\P_\Lambda \bigl( \Omega_\Lambda(\sigma) \ \big|\
\Omega_\Lambda(\mu) \bigr) \leq e^{-\beta \sum_{e: \sigma(e) = 1 }
L(e) -\beta h_{1} \binom{d_\sigma+1}{2}  } e^{\beta h_0
d_{\sigma}}.
\end{equation}
\end{lemma}

The proof of Lemma~\ref{main1} is similar to the proof of Lemma
\ref{main}. The difference in the right hand side between
\reff{ineq1} and \reff{estim} can be explained in the following
way. The set $B$ is nonempty, thus the number of interacted pairs
is at least $\binom{d_\sigma}{2} + d_\sigma = \binom{d_\sigma
+1}{2}$. That provides the energy $h_1 \binom{d_\sigma+1}{2}$.
Removing the star $\sigma$ we can obtain at most $d_\sigma$
isolated points, but not $d_\sigma +1$ as in the \reff{estim},
because now the point $\gamma$ cannot become isolated point.


\subsection{Proof of Theorem \ref{exist} on the existence}
Let $\Gamma$ be weakly homogeneous.
\begin{lemma}\label{mean}
The following inequality
\begin{equation}\label{meanineq}
\E_\Lambda d_{\sigma}(\gamma)\leq \sum_{k=0}^\infty e^{-\beta
h_1\binom{k}{2} + \beta h_0(k +1)}\frac{( T_{\gamma}(\beta)
)^k}{(k-1)!}
\end{equation}
holds for the mean value of the vertex degrees. Where $\E_\Lambda$
is the expectation with respect to the probability $\P_\Lambda$
and $T_{\gamma}(\beta)$ is defined in \reff{hom}.
\end{lemma}

\proof The assertion of the lemma follows from the inequalities
\begin{eqnarray*} \E_\Lambda d_{\sigma}(\gamma) &=&
\sum_{\sigma} d_{\sigma}(\gamma) \P(\Omega_\Lambda(\sigma)) =
\sum_{k=0}^\infty k \sum_{\sigma:\ |E_{\gamma}(\sigma)| = k}
\P(\Omega_\Lambda(\sigma))
\\ &\le& \sum_{k=0}^\infty k e^{-\beta h_1\binom{k}{2} + \beta
h_0(k +1)}\sum_{\sigma:\ |E_{\gamma}(\sigma)| = k} e^{-\beta
\sum_{e\in E_{\gamma}(\sigma)} L(e)} \\ &\le& \sum_{k=0}^\infty k
e^{-\beta h_1\binom{k}{2} + \beta h_0(k +1)}\frac{(
T_{\gamma}(\beta))^k}{k!}.
\end{eqnarray*}
\qed

\begin{rem} \label{rem1} We note that, when $\Lambda$ contains the star, then the
estimation does not depend on $\Lambda.$ Thus it gives an uniform
over $\Lambda$ upper estimation, which holds when $\Lambda
\nearrow \mathbb R^d $
\end{rem}

\vspace{0.5cm}

The theorem \ref{exist} follows now from the finiteness of
$\E_\Lambda d_{\sigma}(\gamma)$ (see \reff{meanineq} and
Remark~\ref{rem1}) \qed

\vspace{0.5cm}

Any sample of Poisson process is weakly homogeneous. It shows the
next

\begin{lemma}\label{poisfin}
Almost all samples $\Gamma$ from Poisson distribution $\pi_{\lambda}$
are weakly homogeneous.
\end{lemma}
\proof Let $\gamma\in\Gamma$. Consider a sequence of rectangles
$U_{n}=[-l_{n},l_{n}]^{d},\;n\geq 0$ centered in $\gamma$ with
size length $l_n.$ We chose $l_{n}=(n+1)^{1/d}$. Then any ring
$W_{n}=U_{n}-U_{n-1},\;n\geq 0$, except $U_{-1}=\emptyset$, has
its volume equal to 1. If $\gamma\in W_{n}$ then for $e=\langle
\gamma,\gamma' \rangle$ the inequality $L(e)\geq l_{n-1}$ holds.
Let $\xi_{n}$ be a number of points from $\Gamma$ located in
$W_{n}$. The variables $\xi_n$ are independent random variables
having Poisson distribution with the parameter $\lambda$ (since
the volume off $W_n$ is equal to 1). Therefore the following
series converges
\begin{equation*}
T_{\gamma}=\sum_{e\in E_{\gamma}}e^{-L(e)}\leq
\sum_{n=0}^{\infty}\xi_{n}e^{-l_{n-1}}
=\sum_{n=0}^{\infty}\xi_{n}e^{-n^{\frac{1}{d}}}<\infty \quad\mbox{
a.s.}
\end{equation*}
The convergence with probability 1 follows from the convergence of
the series of the expectations and the variances of the random
variables $\xi_n e^{-n^{\frac{1}{d}}}$ (Theorem of "two series",
\cite{Shir}). \qed

\subsection{Proof of Theorem \ref{ground} and \ref{uniqgr} on the ground states}

\pr{ground} First we prove the property of the ground states if
there exists at least one. Assume the inverse. Let $\omega$ be a
ground state and there be a vertex $\gamma_1 \in \Gamma \cap
\Lambda$ such that $d_\omega(\gamma_1)\ge 2$. Let $e=[\gamma_1,
\gamma_2]$ be the incident to the vertex $\gamma_1$ in graph
$\omega$, $\omega(e)=1$. Let $\tilde \omega$ be the new
configuration such that $\tilde \omega$ is the same as $\omega$
with the exception that the edge $e$ is now removed: $\tilde
\omega(e)=0$. Then we have
\begin{equation*}
H_\Lambda(\omega)-H_\Lambda(\tilde \omega)=
\begin{cases}
L(e)+ \left(d_{\tilde \omega}(\gamma_1)+d_{\tilde
\omega}(\gamma_2)\right) h_1
& \text{if $d_{\tilde \omega}(\gamma_2)\ge 1$}, \\
L(e)+d_{\tilde \omega}(\gamma_1) h_1 -h_0 & \text{if $d_{\tilde
\omega}(\gamma_2)=0$}.
\end{cases}
\end{equation*}
Since $0<h_0<h_1$, we have
\begin{equation*}
H_\Lambda(\omega)-H_\Lambda(\tilde \omega)>0.
\end{equation*}

There is no edges in a ground with its length $L$ greater than
$2h_0$, since the energy of two monomers is $2h_0<L.$

Proving the existence at least one of the ground states consider a
sequence $(V_n)$ of increasing cubes covering $\R^d=\bigcup_nV_n$.
We build a ground state of the model by a sequence of
reconstructions of an initial configuration. It is reasonable to
take the initial configuration satisfying the property proved
above. For example, we can take the configuration $\omega_0$ with
no edges, that is the configuration of all monomers. Let
$\omega_n$ be a configuration in $V_n$ having the minimal energy
over all configurations in $V_n$. There exists a sequence
$(\omega'_i)$ of configurations which is a subsequence of
$(\omega_n)$, that is $\omega'_i=\omega_{n_i}$, such that there
exists a limit $\lim_{i\to\infty}\omega'_i(e)$ for every $e\in E$.
Moreover, the sequence $(\omega'_i)$ can be chosen such that
$\omega'_j(e)\equiv\const$ for all $j\geq i$ when $e\in E^{V_i}$.
The configuration $\omega'=\bigcup_i\omega'_i$ is one of the
ground states. Indeed, let $\wh{\omega}$ be a local perturbation
of $\omega'$. There exists $V_{i_0}$ such that
$\{e:\:\wh{\omega}(e)\neq\omega'(e)\}\subseteq E^{V_{i_0}}$. For
any $i>i_0$ let $\wh{\omega}_i$ be the configuration equal to the
restriction of $\wh{\omega}$ on $E^{V_{i}}$. The configuration
$\wh{\omega}_i$ is the perturbation of $\omega'_i$ therefore
\begin{equation*}
H_{V_i}(\wh{\omega}_i)-H_{V_i}(\omega'_i) \ge 0.
\end{equation*}
Moreover, the fact, that for any $i$ there is no edges with length
greater that $2h_0$, means that there exists $i_1\ge i_0$ such
that
\begin{equation*}
H(\wh{\omega})-H(\omega') = H_{V_{i_1}}(\wh{\omega}_{i_1}) -
H_{V_{i_1}} (\omega'_{i_1}) \geq 0.
\end{equation*}
That proves that the any local perturbation of $\omega'$ increase
the energy. Thus $\omega'$ is really the ground state. \qed

\pr{uniqgr} The uniqueness follows from  two observations. The
first one is that  there is no edges in the ground state with the
length greater than $2h_0.$ Another observation is that  there
exists a critical intensity $\lambda_c$ such that there is no
boolean percolation with radius $h_0$ for all $\lambda< \lambda_c$
(see \cite{CP}, Theorem 3.3).

Thus, for $\lambda < \lambda_c$ any process configuration $\Gamma$
is an union of finite clusters $\Gamma = \cup_{i = 1} ^\infty
\Gamma_i,\ |\Gamma_i|<\infty,$ and for any $i\ne j$ and for any
$\gamma\in \Gamma_i, \gamma' \in \Gamma_j$ the distance $|\gamma -
\gamma'| > 2h_0.$ There are open edges only inside of the clusters
$\Gamma_i$. Since $\Gamma_i$ are finite there exists a unique
configuration of open edges in every $\Gamma_i$ minimizing the
energy. \qed

\subsection{Proof of Theorem \ref{Poiss} and \ref{core} on non-percolation}
\pr{Poiss} The method of the proof is based on the domination
principle. Namely, we construct a Bernoulli measure $\nu$ on
$\Omega$ which does not percolate and stochastically dominates the
Gibbs measure $\P$. We can apply this method for small rates
$\lambda$ of Poisson measure $\pi_{\lambda}$ and low temperature
of the distribution of Gibbs random graph.

On the set $\Omega$ of the configurations we define the following
Bernoulli measure $\nu$
\begin{equation}\label{b}
\nu( \omega(e) = 1) = \frac{e^{-\beta L(e)}}{e^{-\beta L(e)} +
e^{-2 \beta h_0}}
\end{equation}
independently for any $e\in E$. This measure forms the random -
connected model (see \cite{CP}, ch. 6), which is driven by Poisson
process with the rate $\lambda$ and connected function
\begin{equation} \label{cf} g(x) = \frac{e^{-\beta x}}{e^{-\beta
x} + e^{- 2 \beta h_0 }},
\end{equation}
the probability of two points to be connected on the distance $x$.
The proof of Theorem \ref{Poiss} is a direct application of
Holley's inequality (see \cite{GM}, Theorem 4.8). It is shown in
the next two lemmas.

\begin{lemma}\label{hol}
The following inequality
\begin{equation}
\label{Holly} \P( \omega( e)=1 \mid \omega _{ e} ) \le \nu(
\omega(e)=1 )
\end{equation}
holds for any $\omega _{ e}$, where $\P( \omega( e)=1 \mid \omega
_{ e} ) $ is the Gibbs conditional probability of $(\omega(e)=1)$
given a configuration $\omega_{e}$ out of the edge $e$.
\end{lemma}

{\it Proof of Lemma~\ref{hol}.} Let $e=[\gamma_1,\gamma_2]$. Then
the conditional probability in \reff{Holly} depends on a
configuration on $(E_{\gamma_{1}} \cup E_{\gamma_{2}}) \setminus
\{ e\}:$
 $$ \P(\omega(e)=1 \mid \omega_e)
= \P( \omega (e) =1 \mid \omega_{\gamma_1} \cup \omega_{\gamma_2
}), $$
 where $\omega_{\gamma_{1}}$ and $\omega_{\gamma_{1}}$ are
configurations on $E_{\gamma_{1}} \setminus \{e\}$ and
$E_{\gamma_{2}} \setminus \{e\}$ respectively, and $\cup$ means
the conjugation of the configurations.

Consider three cases:
\begin{enumerate}
  \item $\omega_{\gamma_{1}}=\omega_{\gamma_{2}}\equiv 0$,
  \item $\omega_{\gamma_{1}}\not\equiv 0,\ \omega_{\gamma_{2}}\equiv
  0.$ This case has the symmetrical version $\omega_{\gamma_{1}}\equiv 0,\
  \omega_{\gamma_{2}}\not\equiv 0$.
  \item $\omega_{\gamma_{1}}\not\equiv 0,\ \omega_{\gamma_{2}}\not\equiv
  0$
\end{enumerate}

Case 1. We have
\begin{eqnarray*}
\P(\omega(e) = 1 \mid \omega_{\gamma_{1}}\cup\omega_{\gamma_{2}}
\equiv 0 ) = \frac{ e^{-\beta L(e)}} {e^{-\beta L( e)} + e^{-2
\beta h_0}} = g(L(e)) = \nu( \omega(e) = 1)
\end{eqnarray*}
which means that for the case 1 Holly's inequality holds.

Case 2. Let $E_{\gamma_1} (\omega_{ \gamma_1 })$ ( $E_{\gamma_2}
(\omega_{ \gamma_2 })$ ) be the set of the open edges of
configurations $ \omega_{ \gamma_1}$ $(\omega_{ \gamma_2}).$ Let
$m: = | E_{\gamma_1 } (\omega _{ \gamma_{1}})|.$ Since the edge
$e$ is open then it  interacts with $m$ open edges from
$E_{\gamma_1 } (\omega _{ \gamma_{1}}).$ Then
\begin{eqnarray*}
\P(\omega(e) = 1 \mid \omega_{\gamma_{1}}\cup\omega_{\gamma_{2}} )
&=& \frac{
e^{-\beta L(e) - \beta m h_1}} {e^{-\beta L(e) - \beta m h_1} + e^{-\beta h_0}} \\
&<& \frac{ e^{-\beta L(e)} } {e^{-\beta L(e) } + e^{-\beta h_0}} <
\frac{ e^{-\beta L(e)} } {e^{-\beta L(e) } + e^{- 2\beta h_0}} =
\nu( \omega(e) = 1).
\end{eqnarray*}

Case 3.
\begin{eqnarray*}
\P(\omega(e) = 1 \mid \omega_{\gamma_{1}}\cup\omega_{\gamma_{2}} )
&=& \frac{ e^{-\beta L( e)
- 2 \beta m h_1}} {e^{-\beta L(e) - \beta m h_1} + 1} \\
&<& \frac{ e^{-\beta L(e)} } {e^{-\beta L(e) } + 1} < \frac{
e^{-\beta L(e)} } {e^{-\beta L( e) } + e^{- 2\beta h_0}} = \nu(
\omega(e) = 1)
\end{eqnarray*}
where $m=| E_{\gamma_1} (\omega_{\gamma_{1}})| + | E_{\gamma_2}
(\omega_{\gamma_{2}})|.$ \qed

In the next lemma we find the condition for the non-percolation of
the random - connected model, which dominate the Gibbs
distribution.

\begin{lemma}\label{lemma1}
In the region \reff{region} there is no percolation in the random
- connected model with Poisson rate $\lambda$ and the connection
function (\ref{cf}).
\end{lemma}

\proof The assertion of the lemma is a
consequence of the Theorem 6.1 of \cite{CP}, which claims that a
random - connected model with the connection function (\ref{cf})
does not percolate if
\begin{equation}\label{cond}
\lambda \int_0^\infty g(x)\ed x < 1.
\end{equation}
Note that for any $\beta >0$ and any $h_0>0$ the integral of
$g(x)$ in (\ref{cond}) is finite. We represent the integral in \reff{cond} as
\begin{eqnarray*}
\int g(x)\ed x&=&\int_0^{2h_0}\frac{e^{-\beta x}}{e^{-\beta x} +
e^{- 2 \beta h_0 }}\ed x+ \int_{2h_0}^\infty\frac{e^{-\beta
x}}{e^{-\beta x} + e^{- 2 \beta h_0 }}\ed x \\ & =: & J_1(T) +
J_2(T),
\end{eqnarray*}
where $T=1/\beta.$ The first integral $J_1(T)$ on the right side
of the above equality is increasing and tends to $2h_0$ as
$\beta\rightarrow \infty$. The second integral tends to $0$ as
$\beta\rightarrow \infty$. Choosing $\lambda$ such that
$$ \lambda < \frac{1}{2h_0 + J_2(T)} \le \frac{1}{\int g(x)dx} $$
 we obtain the claim  of the lemma. \qed

\vspace{.5cm}

\pr{Poiss} By Holley inequality, Lemma \ref{hol} implies that the
Gibbs measure on $\Omega$ is dominated by the product measure,
Lemma \ref{lemma1} implies that the product measure does not
percolate under the conditions of the theorem.  \qed

\subsubsection*{The Cluster Branching Process}

The proof of Theorem \ref{core} is based on the construction of a
non-homogeneous { \emph{cluster branching}} process of the edges.

An informal description of the cluster branching process is the
following. Let $B$ be some connected set of open edges which forms
a cluster and let $V$ be the set of vertices in the cluster.
Consider the pair $(V,B)$ as a connected graph. The graph distance
$\rho$ between two vertices is the number of edges in a shortest
path connecting them. Fix a vertex $\gamma_0$ in $V$. For any
$n\in \mathbb N$ a sphere with radius $n$ and center $\gamma_0$ is
$$
V^{(n)} = \{ \gamma\in V:\ \rho(\gamma, \gamma_0) = n \},
\mbox{ where }V^{(0)} = \{\gamma_0\}.
$$
The sequence $\{V^{(i)}\,:\,i=1,2,\dots\}$ is a partition of $V =
\cup_{i=0}^\infty V^{(i)}$. Then $B=\cup_{i=1}^\infty B^{(i)}, $
where
$$
B^{(n)} = \{ e=[w, v] \in B:\ w\in V^{(n-1)}\mbox{ and } v\in
V^{(n-1)} \cup V^{(n)} \}.
$$
We interpret the set $B^{(n)}$ as $n$-th offspring generation of
the ancestor set $V^{(n-1)}$. The set $B^{(n)}$ is a set of 'plant
branches' growing from a set of 'buds' $V^{(n-1)}$. We think
$B^{(n)}$ as the state of a branching cluster process at ``time''
$n$.

This construction leads to an ambiguity since the edge $[w, v]\in
B^{(n)}$ can be the offspring of two ancestors $v$ and $w$ if
$v,w\in V^{(n-1)}$. This problem can be solved by introducing an
order along which the embranchment is controlled. The order of the
branching induces a dependence of the offsprings. Another
peculiarity of the branching cluster process is interactions of
the offsprings having different ancestors. These properties differ
the branching cluster process from the standard branching
processes.

The formal definition of the branching cluster process can be made
in the following way.

{\it Construction of Cluster Branching Process}. Recall that
$E_\gamma$ is the set of all edges incident with the point
$\gamma\in \Gamma$. As before we denote $\Sigma_{\gamma} =
\{0,1\}^ {E_\gamma}$ and $ \Sigma_{\gamma, D} = \{0,1\}^ {E_\gamma
\setminus D}$ the set of all configurations on $E_\gamma$ and
$E_\gamma \setminus D$ correspondingly, where $D$ is some set of
the edges.

The path of the cluster branching process is a sequence of triples
$(B^{(n)},V^{(n)}, S_n)$. The distribution of the cluster
branching process is denoted by $\mathbb P$. The precise
definition is the following. Let $\gamma_0 \in \Gamma$ be the
starting point of a branching process path.

\begin{enumerate}
\item[] {\it Initial stage.} $B^{(0)} := \varnothing,$ ${V}^{(0)}
:= \{ \gamma_0 \}$ and $S_0:= \varnothing.$

\item[] {\it First stage.} Let us choose some set of edges
$B^{(1)}\subseteq E_{\gamma_0}$ which are the offsprings of
$\gamma_0$. With help of $B^{(1)}$ we construct the next objects
\begin{enumerate}
  \item[] $V^{(1)} := \{\gamma:\: [\gamma,
\gamma_{0}] \in B^{(1)}\}$,
   \item[] $S_1:=E_{\gamma_{0}}$.
\end{enumerate}
In order to define the offspring probability ${\mathbb{P}}(B^{(1)}
)$ of the ancestor $\gamma_0$ we introduce the star configuration
 $$\sigma_{\gamma_0}(e)=\begin{cases}
 1,&\mbox{ if }e\in B^{(1)},\\
 0,&\mbox{ if }e\in S_{1}\setminus B^{(1)}
\end{cases}
$$
and
$$
\mu_{1}(e)=\sigma_{\gamma_0}(e).
$$

The path of the one step embranchment is $B_1=B^{(1)}$.

 Then
 $$ {\mathbb{P}}(B^{(1)} ) := \P(\Omega(\sigma_{\gamma_0})),$$
 where $\Omega(\sigma_{\gamma_0})$
is the set of all configurations of $\Omega$ such that its
projection on $E_{\gamma_0}$ coincide with the star-configuration
$\sigma_{\gamma_0}.$ It follows from Theorem~\ref{exist} that  the
number of the offsprings from one point is finite.

\item[] {\it Second stage.} Having $B^{(1)}, {V}^{(1)}, S_1$ we
construct the next generation. Namely, we shall define the objects
$B^{(2)}, {V}^{(2)}, S_2$. We shall do it successively according
to an order in $V^{(1)}$. The order is arbitrary. We need it to
avoid the ambiguity in the definition of ancestors of an offspring
$e=[w,v]$ when $w,v\in {V}^{(1)}$. Let $k_{1}=|V^{(1)}|$. Suppose
that the points in $V^{(1)}$ are enumerated in some way, $V^{(1)}
= \{\gamma^{(1)}_{1}, ..., \gamma^{(1)} _{k_{1}}\}.$ We construct
successively $B^{(2)}_i, {V}^{(2)}_i, S_{2,i},\ i=1,\dots, k_1.$
Let us begin with the first point $\gamma_1^{(1)}.$ Let
$B^{(2)}_1$ be a subset of $ E_{\gamma_1^{(1)}} \setminus S_{1}$
which is a offspring set of $\gamma_1^{(1)}.$ Then
\begin{enumerate}
  \item[]
${V}^{(2)}_1 := \{\gamma:\: [\gamma, \gamma_{1}^{(1)}] \in
B^{(2)}_1\}$;
   \item[] $S_{2,1} := S_1 \cup E_{\gamma_1^{(1)}}$.
\end{enumerate}

Since the set $B^{(2)}_1$ is from
$E_{\gamma_1^{(1)}}\setminus S_1$ the initial point $\gamma_0$
can not belong to  ${V}^{(2)}_1$. However the points from
$V^{(1)}$ may belong to ${V}^{(2)}_{1}$.

In order to define the offspring probability we introduce two configurations:
$$\sigma_{\gamma_1^{(1)}}(e)=\begin{cases}
 1,&\mbox{ if }e\in B^{(2)}_1,\\
 0,&\mbox{ if }e\in E_{\gamma_1^{(1)}} \setminus B^{(2)}_1
\end{cases}$$
and
$$\mu_{2,1}(e)=\begin{cases}
 1,&\mbox{ if }e\in B_{2,1},\\
 0,&\mbox{ if }e\in S_{2,1}\setminus B_{2,1}
\end{cases}
$$
where the path $ B_{2,1} = B^{(1)}\cup B^{(2)}_1$. We have
described two steps of the process: branching from $\gamma_0$ and
from $\gamma_1^{(1)}$.

Further  the upper index denotes the number of a stage and the
lower index if it single denotes the number of a step in the
stage. Double lower indices contain both the step and the stage.

The conditional probability of the offsprings ${B}^{(2)}_1$ of the
ancestor $\gamma_1^{(1)}$ given the {\it environment} $B_{1}$ is
\begin{equation} \label{offsp1}
\mathbb{P} (B^{(2)}_1 \mid B_{1}) := \P(
\Omega(\sigma_{\gamma_1^{(1)} } ) \mid \Omega( \mu_{1}) ).
\end{equation}

Assume we have constructed $B^{(2)}_i,V^{(2)}_{i}, S_{2,i}$ and
also we have $B_{2,i}, \mu_{2,i}, i=1,...,m$, where $m<k_{1}$.
Doing the next branching of the point $\gamma_{m+1}^{(1)} $ choose
some set $B_{m+1}^{(2)}$ from the set $E_{\gamma_{m+1}^{(1)}}
\setminus S_{2,m}$ which means the offspring set of
$\gamma_{m+1}^{(1)}$. Then
\begin{enumerate}
\item[] ${V}^{(2)}_{m+1}= \{\gamma:\: [\gamma, \gamma_{m+1}^{(1)}]
\in B^{(2)}_{m+1}\}$; \item[] $S_{2,m+1} := S_{2,m} \cup
E_{\gamma_{m+1}^{(1)}}$.
\end{enumerate}
Now we obtain the path $B_{2,m+1} = B_{2,m} \cup B_{m+1} ^{(2)}.$

To define the offspring probability introduce the configuration
$$\sigma_{\gamma_{m+1}^{(2)}}(e)=\begin{cases}
 1,&\mbox{ if }e\in B^{(2)}_{m+1},\\
 0,&\mbox{ if }e\in E_{\gamma_{m+1}^{(1)}} \setminus S_{2,m+1}
\end{cases}$$
and configuration
$$
\mu_{2,m+1} := \mu_{2,m} \vee \sigma_{\gamma_{m+1}^{(1)} }.
$$
We use the sign $\vee$ to notate the concatenation of two configurations defined on non-intersected sets.

The conditional  probability of offsprings $B^{(2)}_{m+1}$   of
the ancestor $\gamma_{m+1}^{(1)}$ given $B_{2,m}$ is
\begin{equation} \label{offsp}
\mathbb{P} (B^{(2)}_{m+1} \mid B_{2,m}) := \P(
\Omega(\sigma_{\gamma_{m+1}^{(1)} }) \mid  \Omega( \mu_{2,m} ) )
\end{equation}
Having done the construction for $i=1,\dots, k_1$
we obtain
\begin{enumerate}
 \item[] $B^{(2)} := \cup_{i=1} ^{k_1} B^{(2)}_i$;
 \item[] $V^{(2)} := \cup _{i=1}^{k_1} {V}^{(2)}_{i} \setminus V^{(1)} = \{ \gamma^{(2)} _1, \dots, \gamma ^{(2)} _{k_2} \}$;
 \item[] $B_2:= B_{2,k_1} \mbox{ and } \mu_{2} := \mu_{2,k_1}$;
  \item[] $S_{2}:=S_{2,k_1}.$
\end{enumerate}
Remark that the set $\cup _{i=1}^{k_1} {V}^{(2)}_{i}$ can include
points from $V^{(1)}$. The points from the set  $\left( \cup
_{i=1}^{k_1} {V}^{(2)}_{i} \right) \cap  V^{(1)}$ can not have
offsprings. Therefore they are excluded from the next branching
generation.

\item[] {\it $(n+1)$th stage.} Assume we have constructed
$$B^{(n)},{V}^{(n)}=\{\gamma^{(n)}_{1},...,\gamma^{(n)}_{k_{n}}\},
S_n   \mbox{ and } B_n, \mu_n.$$ Then the next generation $B^{(n
+1)} = \cup_{i=1}^{k_n} B^{(n+1)}_{i}$ is constructed  in the same
way as in the second stage with objects $V^{(n+1)}_i, S_{n+1,i},
B_{n+1,i}, \mu_{n+1,i}$ and $\sigma_{\gamma_{m}^{(n+1)} }$. The
offspring probabilities are defined in the same way
\begin{equation} \label{offsp-n}
\mathbb{P} (B^{(n+1)}_{m+1} \mid B_{n+1,m}) := \P(
\Omega(\sigma_{\gamma_{m+1}^{(n)} }) \mid  \Omega( \mu_{n+1,m} )
).
\end{equation}

\end{enumerate}
It completes the construction of the cluster branching process.

\vspace{1cm}

We show next that the cluster branching processes posses the main
feature of the usual branching processes, namely, if the
expectation of the offspring number of one ancestor is less than 1
then the processes extinct.
\begin{lemma}
Assume that there exists $\varepsilon>0$ such that for $n>1$
either
\begin{equation*}
\mathbb E (|B^{(n)}_i| \mid B_{n,i-1}) \leq 1-\varepsilon \mbox{ when } i>1
\end{equation*}
or
\begin{equation*}
\mathbb E(|B^{(n)}_i|\mid S_{n-1})\leq 1-\varepsilon \mbox{ when }
i=1
\end{equation*}
then
\begin{equation*}
\mathbb E(|B|)<\infty,
\end{equation*}
where $B=\cup_{n=1}^\infty B^{(n)}$
\end{lemma}

\Proof follows from the following equalities
\begin{eqnarray*}
{\mathbb E}\bigl[ |B^{(n)}| \bigr] &=& {\mathbb E}\bigl[
\sum_{k=1}^{\infty}I_{\{|B^{(n-1)}|=k\}}|B^{(n)}| \bigr]
\\
&=& \sum_{k=1}^{\infty} {\mathbb E} \Bigl[ I_{\{|B^{(n-1)}|=k\}}
\sum_{i=1}^{k} |B^{(n)}_i| \Bigr]
\\
&=& \sum_{k=1}^{\infty} \sum_{i=1}^{k} {\mathbb E}  \Bigl[ \mathbb E \Bigl[
I_{\{|B^{(n-1)}|=k\}} |B^{(n)}_i| \, \Bigl|\Bigr.\, B_{n,i-1}  \Bigr] \Bigr] \\
&=& \sum_{k=1}^{\infty} \sum_{i=1}^{k} {\mathbb E}  \Bigl[
I_{\{|B^{(n-1)}|=k\}} \mathbb E
\Bigl[ |B^{(n)}_i| \, \Bigl|\Bigr.\, B_{n,i-1} \Bigr] \Bigr].
\end{eqnarray*}
We used the measurability of the event $(|B^{(n-1)}|=k)$ with
respect to the $\sigma$-algebra  generated by $B_{n,i-1}$. We
adopt the above that $B_{n,0}=B_{n-1}$.

Next we obtain
\begin{eqnarray*}
{\mathbb E}\bigl[ |B^{(n)}| \bigr] &\le& (1-\eps)\sum_{k=1}^{\infty}
k{\mathbb E}  \Bigl[ I_{\{|B^{(n-1)}|=k\}}
\Bigr] \\
&=& (1-\epsilon) {\mathbb E} [ |B^{(n-1)}| ]
\end{eqnarray*}
Thus
$$
{\mathbb E}\bigl[ | B | \bigr] \leq \mathbb E[ |B^{(1)}|
]/\epsilon.
$$ \qed

The definition of the cluster branching process is done in a way
such that any maximal component of any configuration $\omega$ can
be \emph{obtained} as a cluster process path. It means the
following.

Let $\gamma_{0}\in \Gamma$, and $\omega_0$ be some configuration
from $\Omega.$ Let $C_{\gamma_{0}}(\omega)$ be the maximal
connected component of open edges of $\omega_0$ containing
$\gamma_{0}.$ We construct a cluster branching process along
$\omega_{0}$ where $\gamma_{0}$ is the initial point of the
cluster path. The only freedom in the cluster process path
deriving is in the choice of the offsprings. Doing the coupling
with chosen configuration $\omega_0$ we define
$\sigma_{\gamma_i^{(n)}}$ as a projection of $\omega_0$ on the set
of edges $E_{\gamma_i^{(n)}} \setminus S_{n,i-1}.$ Then the
probability 
to have a finite connected component $C_{\gamma_{0}}(\omega_0)$
can be obtained as the probabilities \reff{offsp} of the branching
cluster process path made along $\omega_0$. As a consequence the
following equality holds: for any point $\gamma_0 \in
\Gamma$
\begin{equation}\label{fp} \P( C_{\gamma_0} \ \mbox{\rm is finite})
= \mathbb P ( B^{(n)}\ \mbox{\rm not survives})
\end{equation}

The following lemma finish the prove of the theorem.
\begin{lemma}\label{l2} Let $\Gamma$ be strongly homogeneous
(see \reff{hom1}). Then for any small $\epsilon>0$ there exists
$\beta_0=\beta_0(\epsilon)$ such that for all $\beta > \beta_0$
\begin{equation}\label{meanbr}
{\mathbb E} \bigl[ |B^{(n)}_i| \ \bigl|\bigr.\ B_{n,i-1} \bigr] <
1-\epsilon.
\end{equation}
uniformly over $i,\ n>1$. Here $B_{n,0}=B_{n-1}$.
\end{lemma}
 \proof Let $\gamma_i^{(n-1)}$ be the branching point of which
 offsprings are
$B^{(n)}_i$.  Let the previous path be $B_{n,i-1}$. It follows
from (\ref{ineq1}) that:
\begin{equation}\label{est5}
\mathbb P \Big(|B_{i}^{(n)}|=m \ \Big|\  B_{n,i-1} \Big) \le
e^{-\beta \sum_{e\in B_{i}^{(n)}} L(e) - \beta \binom{m+1}{2}}
e^{\beta h_0 m }.
\end{equation}
Thus,
\begin{eqnarray}
&& {\mathbb E}\Big(|B_{i}^{(n)}| \ \Big|\  B_{n,i-1} \Big) = \sum _{m=1}
^\infty m \sum_{ B_{i}^{(n)} :\ |B_{i}^{(n)}|=m } \mathbb P
\Big(|B_{i}^{(n)}|=m \ \Big|\  B_{n,i-1}\Big) \nonumber
\\ && {} \le \sum _{m=1} ^\infty m \exp \Bigl\{ - \beta h_1
\binom{m+1}2  + \beta h_0 m \Bigr\} \times \nonumber
\\ && {} \ \ \ \  \times \sum_{ B_{i}^{(n)} :\ |B_{i}^{(n)}|=m }
\exp \Bigl\{ - \beta \sum_{e\in B_{i}^{(n)} } L(e) \Bigr\} \nonumber \\
&& {} < \sum _{m=1} ^\infty m \exp \Bigl\{ - \beta h_1 \binom{ m +
1}{2} + \beta h_0 m \Bigr\} \frac{( T_{\gamma_i^{(n+1)}} (\beta)
)^m} {m!}, \label{est6}
\end{eqnarray}
where
 $$ T_{\gamma}(\beta) = \sum_{e \in E_{\gamma}}
 e^{-\beta L(e) }. $$
Since $ T_{\gamma}(\beta)$ are uniformly bounded over $\gamma\in
\Gamma$ the choice of large enough  $\beta$ leads to
\reff{meanbr}. \qed

\section{Conclusions}

\begin{enumerate} \item Since the ground state of Gibbs Random
Graph do not percolate the theorems about the non-percolation show
a kind of "stability" of the ground states.

\item Condition of the existence of an infinite cluster is an open
problem. \end{enumerate}

\section{Acknowledgements}

The work of E.P. was partly supported by CNPq grants 300576/92-7
and 662177/96-7, (PRONEX) and FAPESP grant 99/11962-9, RFBR grants
07-01-92216 and 08-01-00105.

The work of V.S. was partly supported by  FAPESP grant 99/11962-9
and CNPq  grant 306029/2003-0.

The work of A.Ya. was partly supported by E26-170.008-2008
(PRONEX),  "Edital Universal 2006" grant 471925/2006-3 and
306092/2007-7 (CNPq).


\begin{thebibliography}{19}

\bibitem{BB} Bela Bollobas. Random graphs. {\it
Cambridge studies in advanced mathematics.} Second Edition.
Cambrige University Press, 2001.

\bibitem{RD} Rick Durrett. Random graph Dynamics. Cambrige University
Press, October 2006.

\bibitem{CP} Ronald Meester and Rahul Roy. Continuum Percolation. {\it
Series: Cambridge Tracts in Mathematics (No. 119).} Cambrige
University Press, 1996.

\bibitem{GM} Hans-Otto Georgii, Olle Hðaggstrðom, Christian Maes. The random
geometry of equilibrium phases. arXiv:math.PR/9905031 v1 5 May
1999.

\bibitem{Shir}A.N. Shiryaev. Probability. Springer 1996.

\end{thebibliography}
\end{document}